\documentclass{amsart}
\usepackage{graphicx} 
\usepackage{amsfonts}
\usepackage{amsmath,amssymb, amscd,amsbsy, amsthm, enumerate}
\usepackage{mathtools}
\usepackage[unicode]{hyperref}
\usepackage{tikz, xcolor, fancyhdr}
\usepackage{emptypage}
\usepackage[final]{changes}
\newtheorem{theorem}{Theorem}[section]
\newtheorem{lemma}[theorem]{Lemma}

\theoremstyle{definition}

\newtheorem{remark}[theorem]{Remark}

\numberwithin{equation}{section}

\def\Q{{\mathbb Q}}
\def\Z{{\mathbb Z}}
\def\N{{\mathbb N}}
\def\R{{\mathbb R}}

\def\C{\mathbb C}
\def\H{\mathcal H}

\def\O{\mathcal O}
\def\SL2C{\text{SL}_2(\mathbb{C})}
\def\SL2OK{\text{SL}_2(\mathcal{O}_K)}
\def\GL2C{\text{GL}_2(\mathbb{C})}

\title{Bianchi Modular Forms and the Rationality of Periods}

\author[G. Anderson]{Gradin Anderson}\address{Brigham Young University, Provo, UT 84602}\email{gradinmanderson@gmail.com}

\author[P. Harrigan]{Peter Harrigan}\address{Hillsdale College, MI}\email{pharrigan@hillsdale.edu}

\author[L. Hoback]{Louisa Hoback}\address{University of Michigan-Dearborn,  Dearborn, MI 48128}\email{lhoback@umich.edu}

\author[M. Pugh]{McKayah Pugh}\address{University of Northern Colorado, CO}\email{pugh7913@bears.unco.edu}

\author[T.A. Wong]{Tian An Wong}\address{University of Michigan-Dearborn, Dearborn, MI 48128}\email{tiananw@umich.edu}


\begin{document}

\begin{abstract}
Using an explicit Eichler-Shimura-Harder isomorphism, we establish the analogue of Manin's rationality theorem for Bianchi periods and hence special values of $L$-functions of Bianchi cusp forms. This gives a new short proof of a result of Hida in the case of Euclidean imaginary quadratic fields. \added{In particular,} we give an explicit proof using the space of Bianchi period polynomials constructed by Karabulut and describe the action of Hecke operators. \end{abstract}
\maketitle

\section{Introduction}

\subsection{Bianchi period polynomials}

An important method of studying classical modular forms is through their period polynomials.  Let $k\geq0$ be an integer and let $f$ be a cusp form of weight $2k+2$ on SL$_2(\mathbb Z)$. Then the period polynomial of $f$ is defined to be 
\[
r(f)(X) = \int_0^{i\infty} f(z) (X-z)^{2k}dz  = \sum_{n=0}^{2k} i^{-n+1}\begin{pmatrix}2k\\ n\end{pmatrix}r_n(f)X^{2k-n},
\]
and
\[
r_n(f)=\int_0^\infty f(it)t^ndt = n!(2\pi)^{-n-1}L(f,n+1),\qquad 0\le n \le 2k,
\]
is called the $n$-th period of $f$. The map $r(f)$ explicitly realizes the Eichler-Shimura isomorphism, identifying the space of cusp forms with the space of period polynomials. The rationality of periods provides a rich connection between modular forms and arithmetic, e.g. \cite{KOHNENZAGIER}, and remains an active area of study.

Now consider $K=\Q(\sqrt{-D})$, a imaginary quadratic number field  and $\mathcal O_K$ its ring of integers. From an analytic point of view, Karabulut recently constructed a space of Bianchi period polynomials over Euclidean $K$, relating it to the cuspidal cohomology of $\mathrm{SL}_2(\mathcal O_K)$ \cite{KARABULUT}. Combes further defined a Hecke action on period polynomials by relating it to the space of Bianchi modular symbols \cite{COMBES}. On the other hand, the algebraic theory of Bianchi modular forms have also seen new developments, such as  \cite{WILLIAMS,BARRERA,PALACIOS} and many others extending earlier works of Hida \cite{HIDA} and Ghate \cite{GHATE} into the $p$-adic setting for general $K$.

\subsection{Main results} In this paper, we build on these recent advances to study the rationality of periods of Bianchi modular forms.  Let $F$ be a Bianchi cusp form of parallel weight $(k,k)$ on SL$_2(\mathcal O_K)$, we construct the Bianchi period polynomial realizing the Eichler-Shimura-Harder isomorphism in Theorem \ref{thm:ESH},
\[
r(F)(X,Y,\overline{X},\overline{Y})= \int_{0}^{\infty}\omega_F=\sum_{p,q=0}^k \binom{k}{p}\binom{k}{q}r_{p,q}(F)X^{k-p}Y^{p}\overline{X}^{k-q}\overline{Y}^{q}, 
\]
where  $\omega_F$ is the differential 1-form associated to $F$ and $r_{p,q}(F)$ is defined in \eqref{eq:rperiods} as { a period of the Bianchi modular form}
\[
2\binom{2k+2}{k+p-q+1}^{-1}(-1)^{k+q+1}\int_0^\infty t^{p+q}F_{k+p-q+1}(0,t) dt.
\]
The latter encodes the special values of the $L$-function of $F$ by \cite[Theorem 1.8]{WILLIAMS}. When $K$ is Euclidean, we can identify the image of $r(F)$ with the space of period polynomials $\tilde W_{k,k}$, which is a specific quotient of the space $W_{k,k}$ described by \cite{KARABULUT}. In particular, we explicitly construct for the first time the Bianchi period polynomial in $\tilde W_{k,k}$ associated to a given Bianchi cusp form $F$.

As an application of this construction, we establish the rationality of these periods, and hence the rationality of special values of $L$-functions of Bianchi cusp forms. This can be seen as a short proof of a special case of a theorem of Hida \cite[Theorem 8.1]{HIDA}. It is the analogue of Manin's rationality theorem for periods of classical modular forms. 

\begin{theorem}
\label{thm:rationality}
 Let \added{$K$ be a Euclidean imaginary quadratic field}, $F$ be a normalized Hecke eigenform,
and let $K(F)$ be the number field generated by $K$ and the Fourier coefficients of $F$. Then there exists some $\Omega \in \C^{\times}$ such that
    \[\frac{1}{\Omega}r_{p,q}(F)\in K(F)\] for all $0\le p,q\le k$. 
\end{theorem}

\noindent  \replaced{}{In the case where $K$ is Euclidean,} Our proof is explicit, following the classical methods such as laid out in \cite{LANG,COHEN-STROMBERG}, generalized to the Bianchi case using the explicit description of the space of period polynomials $\tilde W_{k,k}$. Namely, we compute the action of Hecke operators on Bianchi periods, then use this to obtain integral formulas for the periods from which we deduce the main theorem. In fact, most of the paper is dedicated to developing the tools in this special case. A short proof for the case of general $K$ was suggested to us by an anonymous referee, the key observation being that it is enough to use the rational structure of the space of polynomials $V_{k,k}$ occurring in the Eichler-Shimura-Harder isomorphism\added{, given in Remark \ref{gen}}.

\subsection{Outline of paper} 

We conclude with a brief outline of the contents of the paper. In Section \ref{section:def}, we introduce notation and the properties of Bianchi modular forms that we will need. 
 In Section \ref{section:integral}, we recall the integration map into the associated first cohomology space and formulate the period map using the Eichler-Shimura-Harder isomorphism in Theorem \ref{thm:ESH}. In Section \ref{section:polynomial}, we express general periods as integral linear combinations of more basic periods using Hurwitz continued fractions, and use this to compute the action of Hecke operators on the periods in Section \ref{section:heck}. Finally, in Section \ref{section:rat}, we apply these results to prove Theorem \ref{thm:rationality}.

\section{Definitions}
\label{section:def}
\subsection{Notation}

Let $K=\Q(\sqrt{-D})$  be an imaginary quadratic number field with class number $1$ 
and let $\O_K$ represent its ring of integers with associated norm $N(\alpha)=\alpha\overline{\alpha}.$
In later definitions, we will also refer to a Euclidean imaginary quadratic number field $K$, where $D$ is $1,2,3,7,$ or 11.
We let 
$$\H_3 = \{(z,t)\in \C\times \R \mid t > 0\}$$
designate hyperbolic 3-space, the unique simply-connected Riemannian manifold of dimension 3 and constant sectional curvature $-1$. 
Additionally, we define the standard action of $\GL2C$ on $\H_3$. Explicitly, for a matrix in $\GL2C$ of determinant $\Delta$, we have 
\[\begin{pmatrix} a & b \\ c & d \end{pmatrix} \cdot (z,t) =  \left( \frac{(az+b)\overline{(cz+d)}+a\bar{c}t^2}{|cz+d|^2+|c|^2t^2},\frac{|\Delta| t}{|cz+d|^2+|c|^2t^2} \right ). \]
Alternatively, we may view elements of $\H_3$ as elements of the Hamiltonians $\mathbb{H}$ via the map $(z,t)\mapsto z+tj$. In this case, for $u\in \mathbb{H}$ the action of $\GL2C$ is given by \[\begin{pmatrix} a & b \\ c & d \end{pmatrix} \cdot u  =  \frac{au+b}{cu+d}. \]
\subsection{Bianchi Modular Forms}

Let $k$ be a non-negative integer. Then denote $V_{2k+2}(\C)$ as the complex vector space of homogeneous polynomials of degree $2k+2$ in variables $X,Y$. 
Define the multiplier system
$$
J(\gamma; (z,t)) = \begin{pmatrix} cz+d & -ct \\ \bar{c}t & \overline{cz+d}\end{pmatrix},\qquad \gamma=\begin{pmatrix} a & b \\ c&d \end{pmatrix} \in \textrm{SL}_2(\mathbb{C}),\ (z,t) \in \H_3.
$$
Given a function $F:\H_3\to V_{2k+2}(\C)$, fixing $(z,t)$ we may represent it as  a polynomial as $F(z,t)=P(X,Y)$ in $V_{2k+2}(\C)$, equipped with the action 
\[
F\left(z,t\right)\left(\gamma\begin{pmatrix} X\\Y\end{pmatrix}\right) = P\left(\gamma\begin{pmatrix} X\\Y\end{pmatrix}\right)= P(aX+bY,cX+dY).
\]
Following \cite{PALACIOS}, for example, given $\gamma \in \mathrm{SL}_2(\mathbb C)$ and a function $F:\H_3 \rightarrow V_{2k+2}(\C)$, we define the slash operator as 
\begin{equation}
\label{eqn:slash}
    F|_{\gamma}(z,t)\begin{pmatrix} X\\Y\end{pmatrix}= |\text{det}(\gamma)|^{-k}F\left(\gamma \cdot (z,t)\right)\left( J\left(\frac{\gamma}{\sqrt{\text{det}(\gamma)}}; (z,t) \right)^{-1}\begin{pmatrix} X \\ Y\end{pmatrix}\right),
\end{equation}
where we note that our $J$ is transposed from that of Palacios. 
Although we may define Bianchi cusp forms for any congruence subgroup of $\SL2OK$, we will only be concerned with those of full level $\SL2OK$. For the remainder of our exposition, we set  $\Gamma = \SL2OK.$ A more general definition is given by \cite[Definition 1.2]{WILLIAMS}.

We say that a function $F:\H_3\rightarrow V_{2k+2}(\C)$ is a \emph{cuspidal Bianchi modular form of weight $(k,k)$ and level $\Gamma$} if it satisfies the following properties:
\begin{enumerate}
    \item[(i)] $F|_\gamma = F$ for every $\gamma \in \Gamma$, 
    \item[(ii)] $\Psi F = 0$ and $\Psi' F = 0$, 
    \item[(iii)] $F$ has at worst polynomial growth at each cusp of $\Gamma$.
    \item[(iv)] $\int_{\C/\O_K} F|_\gamma (z,t) dz = 0$ for every $\gamma \in \Gamma$.
\end{enumerate} 
Here $\Psi, \Psi'$ represent the Casimir operators, which are two elements generating the center of the universal enveloping algebra of the Lie algebra associated to the real Lie group $\textrm{PSL}_2(\mathbb{C})$. If we do not impose (iv) in the above definition, then we say $F$  is a Bianchi modular form of weight $(k,k)$ and level $\Gamma$, although we will not be making use of Bianchi modular forms which are not cuspidal. Finally, we denote by $S_{k,k}(\Gamma)$ to be the space of such Bianchi cusp forms.

\subsection{Fourier-Bessel Expansion}
Let $K_n(x)$ be the modified Bessel function that solves the differential equation
\[\frac{d^2K_n}{dx^2}+\frac1x\frac{dK_n}{dx}-\left(1+\frac{n^2}{x^2}\right)K_n=0.\]
It is well established (see \cite{GHATE} and \cite{PALACIOS}) that a cuspidal Bianchi form $F$ has a Fourier-Bessel expansion given by 

\[F(z,t)\begin{pmatrix} X\\Y\end{pmatrix} = \sum_{n=0}^{2k+2}F_n(z,t)X^{2k+2-n}Y^n,\]
where
\[F_n(z,t) := t\binom{2k+2}{n}\sum_{\alpha\in K^\times} c(\alpha \delta) \left(\frac{\alpha}{i|\alpha|}\right)^{k+1-n}K_{n-k-1}(4\pi |\alpha| t)e^{2\pi i (\alpha z + \overline{\alpha z})},\]
$\delta = \sqrt{-D}$, and the Fourier coefficient $c(I)$ is a function on the fractional ideals of $K$ with $c(I)=0$ for $I$ non-integral.



\section{The period integral}
\label{section:integral}\subsection{The Eichler-Shimura-Harder isomorphism}
Let $V_{k,k}=V_{k,k}(\C)$ denote the complex vector space of polynomials in $X,Y,\overline{X},\overline{Y}$ which are homogeneous of degree $k$ in $X,Y$ and homogeneous of degree $k$ in $\overline{X},\overline{Y}$. Denote by $\mathcal V_{k,k}$ the locally constant sheaf on $Y_\Gamma = \Gamma\backslash \mathbb H$ associated to $V_{k,k}$. Let 
$H^1_{\text{cusp}}(Y_\Gamma, \mathcal V_{k,k})$ denote the cuspidal cohomology of $Y_\Gamma$. The Eichler-Shimura-Harder isomorphism \cite{harder} shows that \begin{equation}
\label{ESH}
S_{k,k}(\Gamma) \cong H^1_{\text{cusp}}(Y_\Gamma, \mathcal V_{k,k}),
\end{equation}
while work by Ghate \cite[\S5]{GHATE} 
gives an explicit form of this isomorphism, which we shall describe.

When $K$ is Euclidean, Karabulut \cite[\S5]{KARABULUT} defines the parabolic group cohomology $H^1_{\text{par}}(\Gamma, V_{k,k})\simeq H^1_{\text{cusp}}(Y_\Gamma, \mathcal V_{k,k})$ which is the \replaced{}{the} quotient of parabolic cocycles
on $\Gamma$ by parabolic coboundaries on $\Gamma$, which we present in the following subsection. Karabulut then gives an isomorphism from $H^1_{\text{par}}(\Gamma, V_{k,k})$ to a quotient of an explicit subspace $W_{k,k}\subset V_{k,k}$, which we denote $\tilde W_{k,k}$. We will combine these results to obtain a map
\begin{equation}
\label{PP}
S_{k,k}(\Gamma) \rightarrow \tilde W_{k,k},
\end{equation}
generalizing the Eichler integral for classical modular forms into the space of period polynomials. 


 Williams \cite[\S2]{WILLIAMS} defines elements of the compactly supported cohomology group $H^1_c(Y_\Gamma, \mathcal V_{k,k})$ by identifying it with the space Symb$_\Gamma(V_{k,k})$ of Bianchi modular symbols,  sending $F$ to $\psi_F$, where $\psi_F$ is the map from pairs of cusps in $\mathbb{P}^1(K)$ to $V_{k,k}$ defined by 
\begin{equation}
\label{psiF}
\psi_F(r,s) = \int_r^s\omega_F,
\end{equation}
and $\omega_F$ is a $V_{k,k}$-valued differential form on $\H_3$ which is defined explicitly. Note that the latter integral gives a perfect pairing between the relevant first homology and cohomology groups. (Williams actually considers spaces $V^*_{k,k}$, the $\C$-dual of $V_{k,k}$, and $V^r_{k,k}$, the space $V_{k,k}$ with $\Gamma$ acting from the right; we will use the latter.) The natural image of $H^1_c(Y_\Gamma, \mathcal V_{k,k})$ in $H^1(Y_\Gamma,\mathcal V_{k,k})\simeq H^1(\Gamma,V_{k,k})$ is often denoted $H_!^1(Y_\Gamma,\mathcal V_{k,k})$ \cite[\S1.4]{ash} and contains $H^1_\text{cusp}(Y_\Gamma,\mathcal V_{k,k})$ \cite[\S5]{HIDA} (c.f. \cite[\S1]{Rahm} for a variation). Putting these together, we have 
\begin{equation}
\label{psi}
S_{k,k}(\Gamma) \simeq H^1_{\text{cusp}}(Y_\Gamma, \mathcal V_{k,k}) \simeq H^1_{\text{par}}(\Gamma, V_{k,k}) \simeq \tilde W_{k,k}
\end{equation}
which we shall make explicit. That is, given a Bianchi modular form $F$, we can associate to it a cuspidal cohomology class such that its image in Symb$_\Gamma(V_{k,k})$ is the modular symbol $\psi_F$ above valued in $V_{k,k}$, and its image in $V_{k,k}$ can be identified with the quotient $\tilde W_{k,k}$.

We first give an expression for $\psi_F$ in a form useful to us. Set the binomial coefficient $(\begin{smallmatrix}n\\ m \end{smallmatrix}) = 0$ whenever $m<0$ and $m>n$.

\begin{lemma}
\label{lem:aperint}
    For $\kappa\in K$, { we have}
    \begin{multline*}
{ \psi_F(\kappa,\infty)=}    \sum_{p,q=0}^k \left(\binom{k}{p}\binom{k}{q}\sum_{i=0}^p\sum_{j=0}^q c_{i,j}(\kappa,F)(-1)^{p+q-i-j}\binom{k-i}{p-i}\binom{k-j}{q-j} \kappa^{p-i}\bar{\kappa}^{q-j} \right) \\ \cdot Y^{p}X^{k-p}\overline{Y}^{q}\overline{X}^{k-q}
    \end{multline*}
    where
    \begin{equation}
    \label{cij}
        c_{i,j}(\kappa,F) = 2\binom{2k+2}{k+i-j+1}^{-1}(-1)^{k+j+1}\int_0^\infty t^{i+j}F_{k+i-j+1}(\kappa,t) dt.
    \end{equation}
\end{lemma}

\begin{proof}
This is essentially \cite[Proposition 2.9]{WILLIAMS}, with the difference being that we have $V_{k,k}$ in place of $V^*_{k,k}$.  We use the following definition of a slash operator on $V_{k,k}$.
    For a polynomial $P\in V_{k,k}$ and a matrix $\gamma = \begin{pmatrix}
        a & b \\ 
        c & d
    \end{pmatrix}\in \GL2C$ we define 
    \begin{equation}
        P(X,Y,\overline{X},\overline{Y})|_\gamma = P\left(aX+bY,cX+dY,\overline{a}\overline{X}+\overline{b}\overline{Y}, \overline{c}\overline{X}+\overline{d}\overline{Y}\right).
    \end{equation}
By \cite[Proposition 2.7]{WILLIAMS}, we have $\psi_F(r,s)=\psi_F(\gamma r,\gamma s)|_\gamma$, and in particular
\begin{equation}
    \label{prop:gammainvarience}
\int_{\gamma(0)}^{\gamma(\infty)} \omega_F = \int_0^\infty \omega_F \Big|_{\gamma^{-1}}.
\end{equation}

Then by \cite[Proposition 2.9]{WILLIAMS} we have
\[
\phi_F(\kappa,\infty) = \sum_{p,q=0}^k c_{p,q}(\kappa,F)(\mathcal{Y}-\kappa\mathcal{X})^{k-p}\mathcal{X}^p(\overline{\mathcal{Y}}-\bar{\kappa}\overline{\mathcal{X}})^{k-q}\overline{\mathcal{X}}^q, 
\]
where $\phi_F=\eta \circ \psi_F$ and $\eta$ is the natural isomorphism from $V^r_{k,k}$ to $V_{k,k}^*$ as given in \cite[Proposition 2.6]{WILLIAMS}, with $\mathcal{X},\mathcal{Y},\overline{\mathcal{X}},\overline{\mathcal{Y}}\in V^*_{k,k}$. We note that the above equality is corrected slightly as explained in \cite[Theorem 3.2]{PALACIOS}.

Rewriting the latter summand as 
\[
\left(\sum_{i=0}^p\sum_{j=0}^q c_{i,j}(\kappa,F)(-1)^{p+q-i-j}\binom{k-i}{p-i}\binom{k-j}{q-j} \kappa^{p-i}\bar{\kappa}^{q-j} \right) \mathcal{X}^{p}\mathcal{Y}^{k-p}\overline{\mathcal{X}}^{q}\overline{\mathcal{Y}}^{k-q},
\]
and undoing the isomorphism $\eta$ via the inverse map
\[\mathcal{X}^{p}\mathcal{Y}^{k-p}\overline{\mathcal{X}}^{q}\overline{\mathcal{Y}}^{k-q} \longmapsto \binom{k}{p}\binom{k}{q}X^{k-p}Y^p\overline{X}^{k-q}\overline{Y}^q \]
then gives the desired expression.
\end{proof}

\subsection{Bianchi period polynomials}

For the rest of this section, we assume $K=\Q(\sqrt{-D})$ is Euclidean so that $D=1,2,3,7,11.$ For each value of $D$, define $\omega = i, i\sqrt{2},\frac{-1+i\sqrt{3}}{2},\frac{1+i\sqrt{7}}{2},$ and $\frac{1+i\sqrt{11}}{2}$ respectively. Furthermore, let
\[S = \begin{pmatrix} 0 & -1 \\ 1 & 0 \end{pmatrix},\quad T= \begin{pmatrix} 1 & 1 \\ 0 & 1 \end{pmatrix},\quad T_\omega= \begin{pmatrix} 1 & \omega \\ 0 & 1 \end{pmatrix},\quad U=TS.\]
When $D=1$, let $L=\begin{pmatrix} i & 0 \\ 0 & -i \end{pmatrix},$ when $D=3$, let $L=\begin{pmatrix} \omega^2 & 0 \\ 0 & \omega \end{pmatrix}$, and when $D=11$, let $E=T_\omega^{-1}ST_\omega ST.$ 

Let $C(\Gamma, V_{k,k})$ be the space of 1-cocycles on $\Gamma$ as in \cite[\S5]{KARABULUT}. Then we define the subspace $C_p(\Gamma, V_{k,k})$ of parabolic { cocycles} by 
\[C_p(\Gamma, V_{k,k})=\{f\in C(\Gamma, V_{k,k})| f(T)=f(T_\omega) = f(L)=0\}\]
when $D=1,3$ and 
\[C_p(\Gamma, V_{k,k})=\{f\in C(\Gamma, V_{k,k})| f(T)=f(T_\omega) =0\}\] 
when $D=2,7,11.$

Karabulut \cite[\S5]{KARABULUT} shows that evaluation of a cocycle at $S$ is an isomorphism from the parabolic cocycles $C_p(\Gamma)$ to a subspace $W_{k,k}$ of $V_{k,k}$ defined for each $D$ as the subspace of polynomials $P\in V_{k,k}$ satisfying  the given relations:
\begin{align*}
&D=1, 3: && P|_{(I+S)}=P|_{(I-L)}= P|_{(I+U+U^2)}
    =P|_{(I+T_\omega SL+(T_\omega SL)^2)}= 0,\\
&D=2: && P|_{(I+S)} = P|_{(I+U+U^2)}=P|_{(I+ST_\omega+T_\omega S+T_\omega^{-1}ST_\omega S)} = 0,\\
&D=7: && P|_{(I+S)} = P|_{(I+U+U^2)}=P|_{(T+ST_\omega+T_\omega ST+ST_\omega^{-1}ST_\omega)} = 0,\\
&D=11: && P|_{(I+S)} = P|_{(I+U+U^2)}
    =P|_{(T+ST_\omega+TE+ST_\omega E^{-1}+T_\omega ST+ST_\omega^{-1}ST_\omega)} = 0.
    \end{align*}
More generally, we may define $W_{k,k}$ as the subspace of $V_{k,k}$ defined by the relations amongst the generators of $\Gamma$.

For the isomorphism with $H^1_\text{par}(\Gamma, V_{k,k})$, it is necessary to identify the image of $B_p(\Gamma)$ of parabolic coboundaries, which Karabulut claims has dimension $1$. We will prove the explicit form of this image, which is also mentioned without proof in \cite{COMBES}.
\begin{lemma}
    The image of $B_p(\Gamma)$ in $W_{k,k}$ is spanned by the polynomial $X^k\overline{X}^k-Y^k\overline{Y}^k$.
\end{lemma}
\begin{proof}
Consider an element $f\in B_p(\Gamma)$. \cite{KARABULUT} shows $f$ is a function from $\Gamma$ to $V_{k,k}$ of the form $\gamma \mapsto P|_\gamma - P,$ where $f(T)=f(T_\omega)=0,$ and additionally $f(L)=0$ for the cases where $d=1,3.$ Let $Q(X,\overline{X})=P(X,1,\overline{X},1).$ We will show that $Q$ is constant.

For this to occur, we must in particular have 
\[
0=Q|_T-Q=Q(X+1,\overline{X}+1)-Q(X,\overline{X}),
\]
\[
0=Q|_{T_{\omega}}-Q=Q(X+\omega,\overline{X}+\bar{\omega})-Q(X,\overline{X}).
\]
Thus $Q(0,0)=Q(n,n)=Q((a+b\omega)n,(a+b\overline{\omega})n)$ for all $n\in \N$ and fixed $a,b\in \N$. Thus shows that $Q(Z,Z)$ and $Q((a+b\omega) Z, (a+b\bar{\omega})Z)$ are both constant polynomials in $Z$ since they attain the same value infinitely many times. In particular, they are the same constant as seen when $Z=0$. This shows that 
\[
Q(Z,Z)=Q\left(\frac{a+b\omega}{a+b\bar{\omega}}Z,Z\right),
\]
so for all $z\in \C^\times$, $Q(w,z)$ attains the same value infinitely many times as long $a,b\in \N$ are chosen so that $\frac{a+b\omega}{a+b\bar{\omega}}$ is not a root of unity.

Therefore, $Q(X,\overline{X})$ is constant as a polynomial in $X$. An analogous argument shows it is also constant as a polynomial in $\overline{X}$, so $Q$ is constant. This shows $P$ is of the form $cY^k\overline{Y}^k$, for some $c\in \C$. Then, $f(\gamma)=P|_\gamma-P$ and its image under evaluation by $S$ is thus 
\[
(cY^k\overline{Y}^k)|_S-cY^k\overline{Y}^k = -c(Y^k\overline{Y}^k-X^k\overline{X}^k),
\] 
and so the result follows.
\end{proof}

Hence, the evaluation at $S$ map is an isomorphism from $H^1_\text{par}(\Gamma, V_{k,k})$ to $W_{k,k}$ modulo the space generated by $X^k\overline{X}^k-Y^k\overline{Y}^k$, which we denote by $\tilde W_{k,k}$. We summarize this as follows.

\begin{theorem}
\label{thm:ESH}
Let $K$ be Euclidean. {The map }
\begin{align*}
{F\mapsto \psi_F(S) =} -2\sum_{p,q=0}^k\binom{k}{p}\binom{k}{q}c_{p,q}(0,F)X^{k-p}Y^{p}\overline{X}^{k-q}\overline{Y}^{q}, 
\end{align*}
where $c_{p,q}(0,F)$ is given in \eqref{cij}, is an isomorphism from $S_{k,k}(\Gamma)$ to $\tilde W_{k,k}$.
\end{theorem}

\begin{proof}
We evaluate the Bianchi modular symbol $\psi_F$ in \eqref{psiF} at the pair of cusps $(r,s) = (0,\infty)$ to obtain an element of $V_{k,k}$, whereby
\[
\psi_F(0,\infty)= \int_{0}^{\infty}\omega_F.
\]
Denoting this temporarily by $P$, we then form the 1-cocycle $P|_\gamma - P$ and evaluate at the element $\gamma = S$ to obtain 
\[
P|_S - P =  \int_{S^{-1}(0)}^{S^{-1}(\infty)}\omega_F - \int_{0}^{\infty}\omega_F = -2\int_{0}^{\infty}\omega_F
\]
which is an element of $\tilde W_{k,k}$. {(This realizes the isomorphism \eqref{psi}, and note that the latter integral is referred to as the canonical period polynomial associated to $F$ in \cite[(5.1)]{COMBES})} Then expanding the integral as in Lemma \ref{lem:aperint} with $\kappa=0$ gives the result.
\end{proof}

\begin{remark} Since we are interested in the integrality of this map, we will rescale and omit the factor of $-2$ in the rest of the paper for the convenience; it does not affect the isomorphism in any material way. 
We shall refer to the image of $F$ under this map the \emph{period polynomial of $F$}. We also write $\psi(F) = \psi_F(S)$ when we wish to emphasize the argument $F$.
\end{remark}


\section{Period relations}
\label{section:polynomial}

\subsection{An intermediate expression}

Let $0\le p,q\le k$ be integers. 
Then from Lemma \ref{lem:aperint} we may write 
\[
\psi_F(\kappa,\infty) = \sum_{p,q=0}^k\binom{k}{p} \binom{k}{q} r_{p,q}(\kappa,F) Y^pX^{k-p}\overline{Y}^{q}\overline{X}^{k-q},\] 
where we denote here
\begin{equation}
\label{eq:rperiods}
r_{p,q}(\kappa, F) = \sum_{i=0}^p\sum_{j=0}^q c_{i,j}(\kappa,F)(-1)^{p+q-i-j}\binom{k-i}{p-i}\binom{k-j}{q-j} \kappa^{p-i}\bar{\kappa}^{q-j}
\end{equation}
and $c_{i,j}(\kappa,F)$ is given in \eqref{cij}. We wish to write the coefficients $r_{p,q}(\kappa,F)$ as an $\O_K$-linear combination of $c_{i,j}(0,F).$
Notice that when $\kappa=0$, $r_{p,q}(0,F)=c_{p,q}(0,F).$ When the context is clear, we will denote $r_{p,q}=r_{p,q}(F)=r_{p,q}(0,F) = c_{p,q}(0,F)$ and $r^{\kappa}_{p,q}=r_{p,q}(\kappa,F)$.

We first compute more general period polynomials.

\begin{lemma}
\label{r0F}
Let $\gamma=\begin{pmatrix} a & b \\ c & d \end{pmatrix}\in \Gamma$ such that $\gamma(\infty) =\infty$. Then 
\[\psi_{F}(\gamma(0),\gamma(\infty)) = \sum_{p,q=0}^k\binom{k}{p} \binom{k}{q} s_{p,q}(\gamma,F) Y^pX^{k-p}\overline{Y}^{q}\overline{X}^{k-q},\]
where $s_{p,q}(\gamma, F)$ is equal to
   \[
    \sum_{i,j=0}^{k}(-1)^{p+q-i-j}e_{i,j}^{p,q}r_{i,j}(0,F),
    \]  
    and
    \[e_{i,j}^{p,q} = \sum_{u,v\in\Z}\binom{p}{u}\binom{k-p}{i-u}a^{u}b^{p-u}c^{i-u}d^{k-p-i+u}\binom{q}{v}\binom{k-q}{j-v}\bar{a}^{v}\bar{b}^{q-v}\bar{c}^{j-v}\bar{d}^{k-q-j+v}.\]
In particular, $e_{i,j}^{p,q}\in \mathcal O_K$.
\end{lemma}

\begin{proof}
Using the identity \eqref{prop:gammainvarience}, we obtain
    \[
    \psi_{F}(\gamma(0),\gamma(\infty)) = \int_{\gamma(0)}^{\gamma(\infty)}\omega_{F} =\left(\int_{0}^{\infty}\omega_{F}\right)\bigg|_{\gamma^{-1}}. 
    \]
Expanding the second integral and applying the $\gamma$ action, we have
   \[
   \sum_{p,q=0}^k \binom{k}{p} \binom{k}{q}r_{p,q}(0,F) (-cX+aY)^p(dX-bY)^{k-p}(-\overline{c}\overline{X}+\overline{a}\overline{Y})^q(\overline{d}\overline{X}-\overline{b}\overline{Y})^{k-q}.
   \]
   The summand can also be written by evaluating the first integral as 
   \[
   \sum_{p,q=0}^k\binom{k}{p} \binom{k}{q}r_{p,q}(\gamma(0),F) Y^{p}X^{k-p}\overline{Y}^{q}\overline{X}^{k-q}.
   \]
   Comparing coefficients, the identity follows.
   \end{proof}

\subsection{Continued Fractions}
In order to transform the endpoints of the period integral as desired, we turn to the Hurwitz continued fractions algorithm for Euclidean imaginary quadratic fields, given in \cite{KARABULUT}.
Let $\kappa\in K$, and let $\lfloor{\kappa_{n}}\rceil$ denote the element of $\O_K$ nearest to $\kappa_n$ in the complex plane with respect to Euclidean distance; to break ties, we round down in both the real and imaginary components. Then our continued fraction algorithm is of the form
\[\kappa_{0}=\kappa, \; \beta_{n}=\lfloor{\kappa_{n}}\rceil, \; \kappa_{n+1}=\frac{1}{\kappa_{n}-\beta_{n}},\] with  $\beta_n\in\O_K$, and for any $\kappa \in K$ is guaranteed to eventually terminate at some point where $\kappa_m=\beta_m$ for some $m$. We define the $n$-th convergent of the continued fraction for $\kappa$ as $\mu_n/\nu_n=\left[\beta_{1},\beta_{2},\ldots,\beta_{n}\right]$, where
\begin{equation}
\label{eq:recursion}
   \begin{aligned}
    \mu_{-2}=0,\hspace{1em} \mu_{-1}=1,\hspace{1em}\mu_{n}=\beta_{n}\mu_{n-1}+\mu_{n-2} \\
    \nu_{-2}=1,\hspace{1em} \nu_{-1}=0,\hspace{1em} \nu_{n}=\beta_{n}\nu_{n-1}+\nu_{n-2}
    \end{aligned}
\end{equation}
and $\lim_{n\to\infty} \frac{\mu_{n}}{\nu_{n}}=\kappa$. Hence for $\kappa=\mu_m/\nu_m$, also denoted as $\kappa=\mu/\nu$, we express \\
\[\frac{\mu}{\nu} =\beta_{0}+\frac{1}{\beta_{1}+\frac{1}{\beta_{2}+\frac{1}{\ddots+\frac{1}{\beta_{m}}}}}=\left[\beta_{1},\beta_{2},\ldots,\beta_{m}\right].\] \\
Analogous to the classical algorithm as in \cite{COHEN-STROMBERG}, we find it useful to express the recursion relations of \eqref{eq:recursion} in matrix form, with
\begin{align*}
    g_{0}=\begin{pmatrix}
 	1 & \beta_{0} \\ 0 & 1
 	\end{pmatrix},\quad  g_{n}=\begin{pmatrix}
 	(-1)^{n}\mu_{n-1} & \mu_{n} \\ (-1)^{n}\nu_{n-1} & \nu_{n}
 	\end{pmatrix} ,\quad 
  g_{n}=g_{n-1}\begin{pmatrix} 0 & (-1)^{n-1} \\ (-1)^{n} & \beta_{n}	\end{pmatrix}.
\end{align*}
Note that det$(g_n)=1$, and hence $g_n\in\Gamma$, for all $n$. Because we retain the same recursion matrices as the classical algorithm with $g_n(0)=\mu_n/\nu_n$ and $g_n(\infty)=\mu_{n-1}/\nu_{n-1}=g_{n-1}(0)$, we can immediately conclude
\begin{equation}
\label{eq:frac}
    \sum_{n=0}^m \int_{g_{n}(0)}^{g_{n}(\infty)}\omega_F=\int_{\frac{\mu_m}{\nu_m}}^{\frac{\mu_{-1}}{\nu_{-1}}}\omega_F=\int_{\frac{\mu}{\nu}}^{\infty}\omega_F.
\end{equation}

\subsection{Integral formulas}

We now show that the coefficients $r_{p,q}(\kappa,F)$ can be given as an $\O_K$-linear combination of periods $r_{i,j}(0,F)$. From there, we also show that the $c_{p,q}(\kappa,F)$ coefficients can be given as a linear combination of the $r_{p,q}(\kappa,F)$ coefficients.

\begin{lemma}
\label{lem:rtheta}For $\kappa\in K$, we have
     \[\psi_F(\kappa,\infty) = \sum_{p,q=0}^k\binom{k}{p} \binom{k}{q} r_{p,q}(\kappa,F) Y^pX^{k-p}\overline{Y}^{q}\overline{X}^{k-q},\] where $r_{p,q}(\kappa,F)$ is an $\O_K$-linear combination of $r_{i,j}(0,F)$.
\end{lemma}
\begin{proof} By \eqref{eq:frac}, we can write $\psi_{F}(\kappa,\infty)$ as 
    \begin{align*}
    \int_{\frac{\mu}{\nu}}^{\infty}\omega_F=\sum_{n=0}^m\int_{g_n(0)}^{g_n(\infty)}\omega_{F}  &=\sum_{n=0}^m\sum_{p,q=0}^k \binom{k}{p} \binom{k}{q} s_{p,q}(g_n,F) Y^{p}X^{k-p}\overline{Y}^{q}\overline{X}^{k-q} \\
    &= \sum_{p,q=0}^k \binom{k}{p} \binom{k}{q} r_{p,q}(\kappa,F) Y^pX^{k-p}\overline{Y}^q\overline{X}^{k-q},
    \end{align*}
equating coefficients and using Lemma \ref{r0F} shows $r_{p,q}(\kappa,F)$ is equal to
\begin{align}
\sum_{n=0}^m \sum_{i,j=0}^{k} (-1)^{p+q-i-j} r_{i,j}(0,F) (e_{i,j}^{p,q})_n,
\label{rpqk}
\end{align}
where $(e_{i,j}^{p,q})_n\in\mathcal O_K$ is equal to
\begin{align*}
 \sum_{u,v \in \Z}
 (-1)^{n(i+j)}\binom{p}{u}\mu_{n-1}^u&\mu_n^{p-u}\binom{k-p}{i-u}\nu_{n-1}^{i-u}\nu_n^{k-p-i+u}
\\
 & \binom{q}{v}\overline{\mu_{n-1}}^v\overline{\mu_n}^{q-v}\binom{k-q}{j-v} \overline{\nu_{n-1}}^{j-v}\overline{\nu_n}^{k-q-j+v}. 
 \end{align*}
Recalling from \eqref{eq:recursion} that $\mu_n,\nu_n\in\O_K$ for all $n$, we have that $r_{p,q}(\kappa,F)$ is an $\mathcal O_K$-linear combination of the periods $r_{i,j}(0,F)$.
\end{proof}

We can now give a formula for the coefficients $c_{p,q}(\kappa,F)$ in terms of periods $r_{i,j}(\kappa,F)$ as follows.
    
    \begin{lemma}
    \label{lem:clincomb}
    For $\kappa\in \mathcal O_K$ (alternatively $K$), the coefficients $c_{p,q}(\kappa,F)$ are an $\O_K$-linear ($K$-linear) combination of periods $r_{i,j}(\kappa,F)$, given explicitly by the formula
        \[c_{p,q}(\kappa,F)=\sum_{i=0}^p\sum_{j=0}^q \binom{k-i}{p-i}\binom{k-j}{q-j}\kappa^{p-i}\bar{\kappa}^{q-j}r_{i,j}(\kappa,F).\]
    \end{lemma}
    \begin{proof}
   
   By \eqref{eq:rperiods} we have
    \[ r_{p,q}(\kappa,F) = \sum_{i=0}^p\sum_{j=0}^q c_{i,j}(\kappa,F)\binom{k-i}{p-i}\binom{k-j}{q-j} (-\kappa)^{p-i}(-\bar{\kappa})^{q-j}. \]
    To invert this linear transformation, consider the polynomial
    \begin{align*}
        &\sum_{p,q=0}^k r_{p,q}(\kappa,F) Y^pX^{k-p}\overline{Y}^q\overline{X}^{k-q} \\
        &= \sum_{p,q=0}^k\left(\sum_{i=0}^p\sum_{j=0}^q c_{i,j}(\kappa,F)\binom{k-i}{p-i}\binom{k-j}{q-j} (-\kappa)^{p-i}(-\bar{\kappa})^{q-j}\right)Y^pX^{k-p}\overline{Y}^q\overline{X}^{k-q} \\
        &= \sum_{p,q=0}^k c_{p,q}(\kappa,F)(X-\kappa Y)^{k-p}Y^{p}(\overline{X}-\bar{\kappa}\overline{Y})^{k-q}\overline{Y}^{q}.
    \end{align*}
    Mapping $X \mapsto X+\kappa Y$ and $\overline{X} \mapsto \overline{X}+\bar{\kappa}\overline{Y}$ yields
    \begin{align*}
        &\sum_{p,q=0}^k c_{p,q}(\kappa,F)X^{k-p}Y^p\overline{X}^{k-q}\overline{Y}^q \\
        &= \sum_{p,q=0}^k r_{p,q}(\kappa,F)(X+\kappa Y)^{k-p}Y^{p}(\overline{X}+\bar{\kappa}\overline{Y})^{k-q}\overline{Y}^{q} \\
        &= \sum_{p,q=0}^k \left(\sum_{i=0}^p\sum_{j=0}^q \binom{k-i}{p-i}\binom{k-j}{q-j}\kappa^{p-i}\bar{\kappa}^{q-j}r_{i,j}(\kappa,F)\right)Y^pX^{k-p}\overline{Y}^q\overline{X}^{k-q},
    \end{align*}
    from which equating coefficients gives the desired result.
    \end{proof}

\begin{remark}
The results and methods of this section easily generalize to $K$ of class number one by the generalized Euclidean algorithm due to Whitley \cite[Section 2.7]{WHITLEY_THESIS}, simply noticing that Whitley's algorithm also produces a sequences of matrices in SL$_2(\mathcal O_K)$.
\end{remark}

\section{Hecke operators}
\label{section:heck}
Now let $K$ be Euclidean and let $\mathfrak p = (\pi)$ be a prime ideal of $\mathcal O_K$. Define the Hecke operator $T_\mathfrak p$ acting on Bianchi modular forms (at full level) as 
\[
F|_{T_\mathfrak p}  = |\pi|^{2k} \sum_{b \in (\mathcal O_K/\pi)^\times }F|_{\left(\begin{smallmatrix}1 & b \\ 0 & \pi\end{smallmatrix}\right)} + F|_{\left(\begin{smallmatrix}\pi & 0 \\ 0 & 1\end{smallmatrix}\right)}.
\]
and extend multiplicatively to all of $\mathcal O_K$ (e.g. \cite{WHITLEY,PALACIOS}). Thus for a Hecke operator $T_\mathfrak n$ for a nonzero integral ideal $\mathfrak n  = n\mathcal O_K$, we have
\begin{equation}
    T_{\mathfrak n}F = |\det(\gamma)|^{2k} \sum_{\gamma\in B_n} F|_\gamma,
\end{equation}
where 
\[B_n = \left\{ \begin{pmatrix}
 	a & b \\ 0 & d
 	\end{pmatrix}\mid ad = n, \text{ and } b \bmod{d}\right\} \]
and the slash operator is defined in \eqref{eqn:slash}.

\begin{theorem}
\label{thm:heckelinearity}
Let $n\in \O_K$. Then $r_{p,q}(0,T_{\mathfrak n}F)$ is an $\O_K$-linear combination of periods of $F$ given by 
    \begin{align*}
\sum_{d|n} \left[\sum_{i=0}^p\sum_{j=0}^q\sum_{\substack{b\bmod d \\ (b,d)=1}} \binom{k-i}{p-i}\binom{k-j}{q-j}b^{p-i}d^i\bar{b}^{q-j}\bar{d}^jr_{i,j}(b/d,F)\right]\left[\sum_{e_1e_2 = \tfrac{n}{d}}e_1^p\bar{e_1}^qe_2^{k-p}\bar{e_2}^{k-q}\right],\end{align*}
    where $r_{p,q}(b/d,F)$ is given in \eqref{rpqk}.
\end{theorem}


\begin{proof} 
{ Using Theorem \ref{thm:ESH}, we first compute }
\[
\psi(T_{\mathfrak n}F) = |n|^{2k}\sum_{\gamma\in B_n}\psi(F|_\gamma),
\]
{by expanding the inner summand as }
\[
\psi\left(|n|^{-k}\sum_{m=0}^{2k+2}\left(\sqrt{\frac{a}{d}}\right)^{2k+2-m}\left(\overline{\sqrt{\frac{a}{d}}}\right)^{m}F_m(\gamma(z,t))X^{2k+2-m}Y^m\right).
\]
Then applying the map $\psi$ and simplifying, we have the $K$-linear sum
\begin{align*}
    & |n|^{-k}\sum_{p,q=0}^k \Bigg(2\binom{2k+2}{k+p-q+1}^{-1}\binom{k}{p}\binom{k}{q}(-1)^{k+q+1}  \\  
    & \cdot \int_0^\infty t^{p+q}\left(\left|\frac{a}{d}\right|^{k+1}\left(\frac{a}{d}\cdot \left|\frac{d}{a}\right|\right)^{q-p}F_{k+p-q+1}(\gamma(0,t))\right)\Bigg)Y^pX^{k-p}\overline{Y}^q\overline{X}^{k-q}.
\end{align*}
Then using the change of variables $t=\frac{|n|t}{|d|^{2}}$, so that $\gamma(0,t)= (\tfrac{b}{d},\tfrac{|n|t}{|d|^2})$, we can write 
\begin{multline*}
r_{p,q}(T_{\mathfrak n}F) = |n|^{k} \sum_{\gamma\in B_n} \Bigg(\left|\frac{a}{d}\right|^{k+1}\left(\frac{a}{d}\cdot \left|\frac{d}{a}\right|\right)^{q-p}\\
\cdot\left[2\binom{2k+2}{k+p-q+1}^{-1}(-1)^{k+q+1}\int_0^\infty t^{p+q}F_{k+p-q+1}\left(\frac{b}{d},\frac{|n|t}{|d|^2}\right)dt\right]\Bigg).
\end{multline*}
By Lemma \ref{lem:clincomb},  the expression in brackets  is equal to
\[
\left|\frac{d^2}{n}\right|^{p+q+1} R_{p,q}^{b/d}:=\left|\frac{d^2}{n}\right|^{p+q+1}  \sum_{i=0}^p\sum_{j=0}^q \binom{k-i}{p-i}\binom{k-j}{q-j}\Big(\frac{b}{d}\Big)^{p-i}\overline{\Big(\frac{b}{d}\Big)}^{q-j}r^{b/d}_{i,j}.
\]
We thus have a $K$-linear expression, which we shall show can be rewritten as an $\mathcal O_K$-linear one.

Now let $\chi(z)=\left( \frac{z}{|z|}\right)^{q-p}$ which is totally multiplicative. Using the definition of $B_n$, we get 
\[
r_{p,q}(T_{\mathfrak n}F)=|n|^{2k-p-q} \sum_{d|n}\sum_{b\bmod d}|d|^{2(p+q-k)}\chi\Big(\frac{n}{d^2}\Big)R_{p,q}^{b/d},
\]
where we say that $d|n$ if there exists some $e\in \mathcal O_K$ such that $n=de$. Writing $b\mapsto be$, $d\mapsto d e$ where $\gcd(b,d)=e$, the inner sum becomes
\[
\sum_{e \mid \tfrac{n}{d}}\sum_{\substack{b\bmod d \\ (b,d)=1}}|de|^{2(p+q-k)}\chi\Big(\frac{n}{d^2e^2}\Big)R_{p,q}^{b/d}.
\]
Rearranging, we see $r_{p,q}(0,T_{\mathfrak n}F)$ is equal to
\[
|n|^{2k-p-q}\chi(n) \sum_{d|n}|d|^{2(p+q-k)}\chi(d^{-2})\left[\sum_{e \mid \tfrac{n}{d}}|e|^{2(p+q-k)}\chi(e^{-2}) \sum_{\substack{b\bmod d \\ (b,d)=1}} R_{p,q}^{b/d} \right],
\]
where the expression in brackets can be rewritten as
\[
d^{-p}\bar{d}^{-q}\left[\sum_{e \mid \tfrac{n}{d}}|e|^{2(p+q-k)}\chi(e^{-2})\right] 
\left[\sum_{i=0}^p\sum_{j=0}^q\sum_{\substack{b\bmod d \\ (b,d)=1}} \binom{k-i}{p-i}\binom{k-j}{q-j}b^{p-i}d^i\bar{b}^{q-j}\bar{d}^jr^{b/d}_{i,j}\right].
\]

Notice that the second term
\[t_{p,q}(d)=\sum_{i=0}^p\sum_{j=0}^q\sum_{\substack{b\bmod d \\ (b,d)=1}} \binom{k-i}{p-i}\binom{k-j}{q-j}b^{p-i}d^i\bar{b}^{q-j}\bar{d}^jr^{b/d}_{i,j}\]
is an $\O_K$-linear combination of periods of $F$. We then see that the first term is equal to
\[
\sum_{e \mid \tfrac{n}{d}}|e|^{2(p+q-k)}\chi(e^{-2}) = \sum_{e \mid \tfrac{n}{d}}\left|\frac{n}{de}\right|^{2(p+q-k)}\chi\left(\left(\frac{n}{de}\right)^{-2}\right).
\]
From this, we find
\begin{align*}
    r_{p,q}(T_{\mathfrak n}F) 
    &= |n|^{p+q}\chi(n^{-1}) \sum_{d|n}d^{-p}\bar{d}^{-q}t_{p,q}(d)\sum_{e \mid \tfrac{n}{d}}|e|^{2(k-p-q)}\chi(e^2) \\
    &= n^{p-q}|n|^{2q} \sum_{d|n}d^{-p}\bar{d}^{-q}t_{p,q}(d)\sum_{e \mid \tfrac{n}{d}}|e|^{2(k-2q)}e^{2(q-p)} \\
    &= \sum_{d|n} t_{p,q}(d)\sum_{e \mid \tfrac{n}{d}}n^p\bar{n}^q d^{-p}\bar{d}^{-q}e^{k-2p}\bar{e}^{k-2q}.
\end{align*}    
Since the inner sum of the above is equal to
    \[
\sum_{e \mid \tfrac{n}{d}}\left(\frac{n}{de}\right)^p\overline{\left(\frac{n}{de}\right)}^qe^{k-p}\bar{e}^{k-q} =     \sum_{e_1e_2 = \tfrac{n}{d}}e_1^p\overline{e_1}^qe_2^{k-p}\overline{e_2}^{k-q}
    \]
which is an element of $\O_K$, we find $r_{p,q}(T_{\mathfrak n}F)$ is an $\O_K$ linear combination of $r_{i,j}^{b/d}(F)$ and hence the periods $r_{i,j}(F)$ by Lemma \ref{lem:rtheta}.
\end{proof}

\section{Rationality of periods}
\label{section:rat}

\subsection{Combinatorial lemmas}
We make use of our previous discussion on Hecke operators to prove an analogue of Manin's rationality result for periods of classical modular forms. 
We will prove the main Theorem \ref{thm:rationality} following the strategy of the proof of the classical case \cite[Theorem 11.11.2]{COHEN-STROMBERG}. First, for a cusp form $F$, we define
\[
r(F)=\left(r_{0,0}(F),r_{0,1}(F),r_{0,2}(F), \cdots , r_{k,k}(F)\right)^\mathsf{T},
\]
where we interpret $r(F)$ as a column vector mapping $S_{k,k}(\Gamma)$ to $\C^{(k+1)^2}$. { For any $\mathfrak n = n\mathcal O_K$, note} that since $r_{p,q}(T_{\mathfrak n}F)$ is an $\O_K$-linear combination of $r_{i,j}(F)$, then  there exists some matrix $A(n)$ with entries in $\O_K$ such that $r(T_{\mathfrak n}F)=A(n)r(F)$. Additionally, \replaced{}{for any $\mathfrak n = n\mathcal O_K$} let 
\[
\Tilde{\sigma}_{k}({ n})=\sum_{d\mid n} |d|^k
\]
similar to the standard sum of divisors function.

First we note some observations about $r_{0,0}(T_{\mathfrak n}F)$ in terms of the following lemma.
\begin{lemma}
\label{lem:r00}
    The first component of $A(n)(1,0,\cdots, 0, -1)^\mathsf{T}$ is $\Tilde{\sigma}_{2k+2}(n)$.
\end{lemma}
\begin{proof}
By Theorem \ref{thm:heckelinearity}, we have 
\begin{align*}
    r_{0,0}(T_{\mathfrak n}F) 
    &= \sum_{d|n} \left[\sum_{\substack{b\bmod d \\ (b,d)=1}} r^{b/d}_{0,0}\right]\left[\sum_{e\mid \tfrac{n}{d}}e^k\bar{e}^k\right] = \sum_{d|n} \Tilde{\sigma}_{2k}\left(\frac{n}{d}\right)\sum_{\substack{b\bmod d \\ (b,d)=1}} r^{b/d}_{0,0}.
\end{align*}
Let $b_\ell/d_\ell$ represent the $\ell$th convergent of $b/d$ and let $m$ be the index such that $b_m/d_m=b/d.$ When $p=q=0$, the proof of Lemma \ref{lem:rtheta} gives
\begin{align*}
    r_{0,0}^{b/d} 
    &= \sum_{\ell=0}^m \sum_{i,j=0}^{k}  \left[(-1)^{(n-1)(i+j)}
\binom{k}{i}d_{\ell-1}^{i}d_\ell^{k-i}\binom{k}{j}\overline{d_{\ell-1}}^{j}\overline{d_\ell}^{k-j}\right]r_{i,j}(F); 
\end{align*}
putting these results together gives that $r_{0,0}(T_{\mathfrak n}F)$ is equal to
\[ \sum_{d|n} \Tilde{\sigma}_{2k}\left(\frac{n}{d}\right)\sum_{\substack{b\bmod d \\ (b,d)=1}} \sum_{\ell=0}^m \sum_{i,j=0}^{k}  \left[(-1)^{(n-1)(i+j)}
\binom{k}{i}d_{\ell-1}^{i}d_\ell^{k-i}\binom{k}{j}\overline{d_{\ell-1}}^{j}\overline{d_\ell}^{k-j}\right]r_{i,j}(F).\]

In particular, the coefficient of $r_{0,0}(F)$ within the linear combination for $r_{0,0}(T_{\mathfrak n}F)$ is given by 
\[
\sum_{d|n} \Tilde{\sigma}_{2k}\left(\frac{n}{d}\right)\sum_{\substack{b\bmod d \\ (b,d)=1}} \sum_{\ell=0}^m 
|d_\ell|^{2k},\]
and the coefficient of $r_{k,k}(F)$ within the linear combination for $r_{0,0}(T_{\mathfrak n}F)$ is given by 
\[
\sum_{d|n} \Tilde{\sigma}_{2k}\left(\frac{n}{d}\right)\sum_{\substack{b\bmod d \\ (b,d)=1}} \sum_{\ell=0}^m 
|d_{\ell-1}|^{2k}.\]

Finally, note that since $\psi(F)\in W_{k,k}$, we have $\psi(F)|_{(I+S)}=0$. This implies $r_{0,0}(F)=r_{k,k}(F)$ and so the difference between the coefficients are given by
\begin{align*}
& \sum_{d|n} \Tilde{\sigma}_{2k}\left(\frac{n}{d}\right)\sum_{\substack{b\bmod d \\ (b,d)=1}} \sum_{\ell=0}^m 
(|d_\ell|^{2k}-|d_{\ell-1}|^{2k}) = \sum_{d|n} \Tilde{\sigma}_{2k}\left(\frac{n}{d}\right)\sum_{\substack{b\bmod d \\ (b,d)=1}} |d|^{2k}.
\end{align*}
At this point, let 
\[\Tilde{\phi}(n) = \left|\left\{ m \in \O_K: N(m)<N(n), (m,n)=1 \right\}\right|.\] Note that since $\O_K$ is a Euclidean domain, $\Tilde{\phi}(n)$ is also equal to $|(\O_K/n\O_K)^\times|$, and hence, is multiplicative by the Chinese Remainder Theorem for number fields.
We claim that
    \[\sum_{d|n} \Tilde{\phi}(d) = N(n)= |n|^2.\]
Let $\pi$ denote a prime of $\O_K$, and $e_\pi$ the largest power of $\pi$ dividing $n$. We have
    \[\sum_{d|n} \Tilde{\phi}(d) = \prod_{\pi|n}\sum_{\nu=0}^{e_\pi} \Tilde{\phi}(\pi^\nu).\]
    Since the norm is multiplicative, it suffices to show that $\sum_{\nu=0}^{e_\pi} \Tilde{\phi}(\pi^\nu)=N(\pi^{e_\pi}).$ We then compute $\Tilde{\phi}(\pi^\nu)$ by 
    \begin{align*}
        &\left|\left\{ m \in \O_K: N(m)<N(\pi^\nu) \right\}\right| - \left|\left\{ m \in \O_K: N(m)<N(\pi^\nu), \pi | m \right\}\right| \\
        &= \left|\left\{ m \in \O_K: N(m)<N(\pi^\nu) \right\}\right| - \left|\left\{ m \in \O_K: N(m)<N(\pi^{\nu-1}) \right\}\right| \\
        &= N(\pi^\nu) - N(\pi^{\nu-1}),
    \end{align*}
    and so 
    \[\sum_{\nu=0}^{e_\pi} \Tilde{\phi}(\pi^\nu) = \sum_{\nu=0}^{e_\pi}\left(N(\pi^\nu) - N(\pi^{\nu-1})\right) = N(\pi^{e_\pi}).\]
Using this, we get our desired result as

\begin{align*}
    \sum_{d|n} \Tilde{\sigma}_{2k}\left(\frac{n}{d}\right)\sum_{\substack{b\bmod d \\ (b,d)=1}} |d|^{2k} &= \sum_{d|n} \Tilde{\sigma}_{2k}\left(\frac{n}{d}\right) |d|^{2k} \Tilde{\phi}(d) = \sum_{d|n} \sum_{e\mid \tfrac{n}{d} }|e|^{2k} |d|^{2k} \Tilde{\phi}(d) \\
    &= \sum_{D|n} |D|^{2k} \sum_{d\mid D } \Tilde{\phi}(d) = \sum_{D|n} |D|^{2k+2} = \Tilde{\sigma}_{2k+2}(n).\\ 
\end{align*}
\end{proof}

\subsection{Proof of Theorem \ref{thm:rationality}}

We are now ready to prove the rationality theorem for the Bianchi case. Let $F$ be a normalized Hecke-eigenform. It follows that $T_{\mathfrak n}(F)=c(n\delta)(F)$, with $\delta = \sqrt{-D}$. Thus 
\[ r(F) \in E_{F}= \bigcap_{{ n\in \O_K\backslash\{0\}}} \text{Ker}(A(n)-c(n\delta)\mathit{I}),\]
defined over { $K(F)$}.
We shall identify the algebraic closure of $K$ with a subset of $\C$, so that $r$ maps $S_{k,k}(\Gamma)$ to a subspace $R_k$ of $\C^{(k+1)^2}$, where the relations defining this subspace are completely determined by the action of $\Gamma$ on periods. It follows then that we can identify $R_k$ with $W_{k,k}$, and moreover $r(F) \in E_F \cap W_{k,k}$, a $\C$-vector space. \\

\begin{lemma}
\label{lemvrg}
Let $K$ be Euclidean. If $v\in E_F \cap W_{k,k}$, then there exists a Bianchi cusp form $G\in S_{k,k}(\Gamma)$ such that $v=r(G)$. 
\end{lemma}

\begin{proof}
By Theorem \ref{thm:ESH}, we have a map from $S_{k,k}(\Gamma)$ to the space of period polynomials $\tilde W_{k,k} = W_{k,k}/U$ where $U$ is spanned by $X^k\overline{X}^k-Y^k\overline{Y}^k$. This can be represented by the coefficient vector $(x, 0, \ldots , 0, -x)^\mathsf{T}$ generated by the cocycle $\gamma \mapsto (1|_{\gamma}-1)$. Thus if $v\in E_F \cap W_{k,k}$, then there exists a Bianchi cusp form $G\in S_{k,k}(\Gamma)$ such that 
\[
v=r(G)+(x, 0, \ldots , 0, -x)^\mathsf{T}.
\]
Consider the first component of the product
\[A(n)(x,0,\dots, -x)^\mathsf{T} = A(n)(v-r(G)).\] 
From Lemma \ref{lem:r00}, this is equal to $x\Tilde{\sigma}_{2k+2}(n)$. However, we also see that 
\[A(n)(v-r(G))=A(n)v-A(n)r(G)=a(n)v-r(T_{\mathfrak n}G),\]
where $a(n)$ is the $n$-th Fourier coefficient and $T_{\mathfrak n}$ eigenvalue of $G$, and has the first component 
\[a(n)(r_{0,0}(G)+x)-r_{0,0}(T_{\mathfrak n}G).\]
By the recent proof of the Ramanujan conjecture for cuspidal Bianchi eigenforms \cite{RAMA}, we have that $|a(n)|=O(N(n)^{(k+1)/2})$ (though any weaker nontrivial bound suffices). Furthermore, it is known \cite[Theorem 2.6.2]{garrett} that 
\[
G=\sum_{i=1}^g \lambda_i F_i
\]
for $F_i$ normalized eigenforms.
Therefore, \[r_{0,0}(T_{\mathfrak n}G)=\sum_{i=1}^g \lambda_i r_{0,0}(T_{\mathfrak n}F_i) = \sum_{i=1}^g \lambda_i a_i(n)r_{0,0}(F_i),
\]
where $a_i(n)$ is the Fourier coefficient of $F_i$, which shows that $|r_{0,0}(T_{\mathfrak n}G)| = O(N(n)^{(k+1)/2})$ also. But since $|x\Tilde{\sigma}_{2k+2}(n)| = O(N(n)^{k+1})$, this is a contradiction unless $x=0$. Therefore, $v=r(G)$ for some $G\in S_{k,k}(\Gamma)$.

\end{proof}

We now complete the proof of the main theorem.  First assume that $K$ is Euclidean. If $v\in E_F \cap W_{k,k}$, then by Lemma \ref{lemvrg} there exists a Bianchi cusp form $G\in S_{k,k}(\Gamma)$ such that $v=r(G)$. Since $v\in E_F$, we have by definition that $A(n)v = c(n\delta)v$, and by the definition of $A(n)$ that \[r\left(T_{\mathfrak n}G-c(n\delta)G\right) = 0.\] 
On the other hand, it follows from our realization of the Eichler-Shimura-Harder in Theorem \ref{thm:ESH} and the proof of the previous lemma that $S_{k,k}(\Gamma)$ is isomorphic to a codimension 1-subspace of $R_k$ not containing the vector $(x,0,\dots, -x)^\mathsf{T} $. In particular, $r$ is injective and 
thus\[T_{\mathfrak n}G=c(n\delta)G.\]
It follows by multiplicity one (e.g. \cite[p.432]{HIDA}) that if two eigenfunctions for all $T_{\mathfrak n}$ have the same eigenvalues, then they must be scalar multiples of one another, so $G=\lambda F$ for some $\lambda \in \C $. Thus every element of $E_F \cap W_{k,k}$ has the form $r(\lambda F) = \lambda r(F)$ by the definition of $r$ and the properties of $r_{i,j}$. But from this it follows that $ E_F \cap W_{k,k}$ has dimension one, and hence 
\[
\langle r(F)\rangle_\C = E_F \cap W_{k,k}.
\] 
Note however that $E_F$, $W_{k,k}$, and $T_{\mathfrak n}$ are all defined over $K(F)$. Thus for each $r_{p,q} (F)$ there exists  $\Omega \in \C^{\times}$ such that \[\frac{1}{\Omega}r_{p,q}(F) \in K(F), \forall \text{ } 0 \leq p,q \leq k.\]

{\begin{remark}\label{gen}
Let $K$ be an arbitrary imaginary quadratic field with class number $h$. We briefly outline the idea of the proof of the main Theorem \ref{thm:rationality}. Consider the intersection 
\[
E'_{F}= \bigcap_{(0)\neq \mathfrak n \subset  \O_K} \text{Ker}(T_{\mathfrak n}-c(\mathfrak n)) \subset \text{Symb}_{\mathrm{GL}_2(\hat{\mathcal{O}}_K)}(V_{k,k})
\]
in the space of modular symbols valued in $V_{k,k}$, with notation as in \cite[\S2]{WILLIAMS}. Let $W_{k,k}$ be the subspace of $V_{k,k}$ determined by the Eichler-Shimura-Harder isomorphism, evaluating \eqref{psiF} at the pair of cusps $(r,s) = (0,\infty)$ as before. Again by multiplicity one and Eichler-Shimura-Harder, $E'_{F}\cap W_{k,k}$ is 1-dimensional and spanned by $\psi_F=(\psi_{F^1},\dots,\psi_{F^h})$ where $F^1,\dots,F_h$ are the components of the Bianchi modular form $F$ as in \cite[Proposition 2.7]{WILLIAMS}. As before, the action of $T_{\mathfrak n}$ on modular symbols over $K(F)$ is again defined over $K(F)$, so the intersection 
\[
\text{Symb}_{\mathrm{GL}_2(\hat{\mathcal{O}}_K)}(V_{k,k}(K(F)))\cap E'_F
\]
is again 1-dimensional and spanned by $\psi_F/\Omega$ for some fixed $\Omega\in \C^\times$ (see \cite[Proposition 2.12]{WILLIAMS}). We then obtain 
$\psi_{F^i}(0,\infty)/\Omega \in V_{k,k}(K(F))$ for each $i = 1,\dots,h$.

\end{remark}
}

\subsection*{Acknowledgments}  This research was completed at the REU Site: Mathematical Analysis and Applications at the University of Michigan-Dearborn. The authors were supported by NSF grants  DMS-1950102 and DMS-2243808, and NSA H98230-23. The last author was partially supported by NSF grant DMS-2212924, and thanks Cihan Karabulut for helpful conversations regarding his work.

\bibliographystyle{alpha}
\bibliography{bibFile}

\newcommand{\etalchar}[1]{$^{#1}$}
\begin{thebibliography}{BCG{\etalchar{+}}25}

\bibitem[AS86]{ash}
Avner Ash and Glenn Stevens.
\newblock Cohomology of arithmetic groups and congruences between systems of
  {H}ecke eigenvalues.
\newblock {\em J. Reine Angew. Math.}, 365:192--220, 1986.

\bibitem[BCG{\etalchar{+}}25]{RAMA}
George Boxer, Frank Calegari, Toby Gee, James Newton, and Jack~A. Thorne.
\newblock The {R}amanujan and {S}ato--{T}ate conjectures for {B}ianchi modular
  forms.
\newblock {\em Forum Math. Pi}, 13:Paper No. e10, 65, 2025.

\bibitem[BSW21]{BARRERA}
Daniel Barrera~Salazar and Chris Williams.
\newblock Families of {B}ianchi modular symbols: critical base-change
  {$p$}-adic {$L$}-functions and {$p$}-adic {A}rtin formalism.
\newblock {\em Selecta Math. (N.S.)}, 27(5):Paper No. 82, 45, 2021.

\bibitem[Com24]{COMBES}
Lewis Combes.
\newblock Bianchi period polynomials: {H}ecke action and congruences.
\newblock {\em Res. Number Theory}, 10(2):Paper No. 40, 23, 2024.

\bibitem[CS17]{COHEN-STROMBERG}
Henri Cohen and Fredrik Str{\"o}mberg.
\newblock {\em Modular Forms: A Classical Approach}.
\newblock American Mathematical Society, 2017.

\bibitem[CW94]{WHITLEY}
J.~E. Cremona and E.~Whitley.
\newblock Periods of cusp forms and elliptic curves over imaginary quadratic
  fields.
\newblock {\em Math. Comp.}, 62(205):407--429, 1994.

\bibitem[Gar18]{garrett}
Paul Garrett.
\newblock {\em Modern analysis of automorphic forms by example. {V}ol. 2},
  volume 174 of {\em Cambridge Studies in Advanced Mathematics}.
\newblock Cambridge University Press, Cambridge, 2018.

\bibitem[Gha99]{GHATE}
Eknath Ghate.
\newblock {Critical values of the twisted tensor $L$-function in the imaginary
  quadratic case}.
\newblock {\em Duke Mathematical Journal}, 96(3):595 -- 638, 1999.

\bibitem[Har87]{harder}
G.~Harder.
\newblock Eisenstein cohomology of arithmetic groups. {T}he case {${\rm
  GL}_2$}.
\newblock {\em Invent. Math.}, 89(1):37--118, 1987.

\bibitem[Hid94]{HIDA}
Haruzo Hida.
\newblock On the critical values of {$L$}-functions of {${\rm GL}(2)$} and
  {${\rm GL}(2)\times{\rm GL}(2)$}.
\newblock {\em Duke Math. J.}, 74(2):431--529, 1994.

\bibitem[Kar22]{KARABULUT}
Cihan Karabulut.
\newblock From binary {H}ermitian forms to parabolic cocycles of {E}uclidean
  {B}ianchi groups.
\newblock {\em J. Number Theory}, 236:71--115, 2022.

\bibitem[KZ84]{KOHNENZAGIER}
W.~Kohnen and D.~Zagier.
\newblock Modular forms with rational periods.
\newblock In {\em Modular forms ({D}urham, 1983)}, Ellis Horwood Ser. Math.
  Appl.: Statist. Oper. Res., pages 197--249. Horwood, Chichester, 1984.

\bibitem[Lan95]{LANG}
Serge Lang.
\newblock {\em Introduction to modular forms / Serge Lang}.
\newblock Springer-Verlag Berlin ; New York, 1995.

\bibitem[Pal23]{PALACIOS}
Luis~Santiago Palacios.
\newblock Functional equation of the {$p$}-adic {$L$}-function of bianchi
  modular forms.
\newblock {\em Journal of Number Theory}, 242:725--753, 2023.

\bibitem[R{\c{S}}13]{Rahm}
Alexander~D. Rahm and Mehmet~Haluk {\c{S}}eng{\"{u}}n.
\newblock On level one cuspidal {B}ianchi modular forms.
\newblock {\em LMS J. Comput. Math.}, 16:187--199, 2013.

\bibitem[Whi90]{WHITLEY_THESIS}
Elise Whitley.
\newblock Modular forms and elliptic curves over imaginary quadratic number
  fields.
\newblock 1990.

\bibitem[Wil17]{WILLIAMS}
Chris Williams.
\newblock {$P$}-adic {$L$}-functions of bianchi modular forms.
\newblock {\em Proceedings of the London Mathematical Society},
  114(4):614--656, jan 2017.

\end{thebibliography}

\end{document}